\input amstex
\documentstyle{amsppt}
\magnification=\magstep1

\baselineskip=13pt
\parskip3pt

\define\m1{^{-1}}
\define\ov1{\overline}
\def\gp#1{\langle#1\rangle}

\topmatter
\title
Covering theorems for Artinian rings
\endtitle
\author A.~BORB\'ELY,  V. BOVDI, B. BRINDZA, T.~KRAUSZ
\endauthor
\abstract
The covering properties of Artinian rings which depend on
their additive structure only, are investigated.
\endabstract

\address
       A. Borb\'ely\newline
       Department of Mathematics and Computer Science\newline
       Faculty of Science\newline
       P.O.Box 5969, Safat 13060\newline
       Kuwait
\endaddress
\address
       V. Bovdi\newline
       Institute of Mathematics and Informatics\newline
       Lajos Kossuth University \newline
       H-4010  Debrecen, P.O.Box 12\newline
       Hungary
\endaddress
\email
vbovdi\@math.klte.hu
\endemail
\address
       B. Brindza\newline
       Institute of Mathematics and Informatics\newline
       Lajos Kossuth University \newline
       H-4010  Debrecen, P.O.Box 12\newline
       Hungary
\endaddress
\address
      T. Krausz\newline
      Centre of Informatics and Computing\newline
      Lajos Kossuth University\newline
      H-4010  Debrecen, P.O.Box 58\newline
      Hungary
\endaddress

\thanks
Research was supported by OTKA  No.~T~014279, No.~F~15470.
\endthanks

\endtopmatter

\document

\head {1. Introduction} \endhead

For simplicity, it is convenient to introduce the following notation. A set
$S$ is said to be the proper union of the sets $S_1, \ldots, S_n$ if
$$
\bigcup_{i=1}^nS_i=S\quad \text{and}\quad \bigcup_{i\ne k}S_i\ne S, $$
for all $k=1, \ldots, n$.
Generalizing some earlier results of \cite{1} and  \cite{2}, like a
field is not a proper union of subfields, {\smc \^Ohori} \cite{3} proved
that if a unitary ring $A$ contains a unitary subring $B$ such that
$B/{\frak J}(B)$, where ${\frak J}(B)$ is the Jacobson radical of
$B$, is an infinite (left) Artinian simple ring then $A$ is not a
proper union of rings. As it was remarked by the reviewer of~\cite{3}
(see~\cite{4}) the word ``Artinian" can be deleted by using a
theorem of {\smc Lewin}~\cite{5}.

The purpose of this note is to point out that the covering
properties of Artinian rings depend on their additive structure and
in case of fields the multiplicative structure can be treated as well.

\head {2. Results} \endhead

\proclaim {Theorem 1}
An Artinian ring is not a proper union of additive subgroups if
and only if its additive group is a direct sum of a divisible group
and a finite cyclic group.
\endproclaim

\proclaim {Corollary 1}
An Artinian ring is not a proper union of cosets if
and only if its additive group is  a divisible group.
\endproclaim

\proclaim {Theorem 2}
A ring with minimal condition for principal left ideals is not a proper union
of cosets if and only if its additive group is a direct sum of a divisible
group and a torsion   group  which has no subgroup of finite index.
\endproclaim

\proclaim {Corollary 2}
A ring with minimal condition for principal left ideals is not a proper union
of additive subgroups if and only if its additive group is a direct sum of a
divisible group, and a torsion  group such that every finite factorgroup of it
is cyclic.
\endproclaim

\proclaim {Theorem 3}
Let $R$ be an infinite skew  field and $\{H_{1},H_{2},\ldots,H_{t}\}$ be a
family of its  proper subskew fields. Then
\item {\rm (i)}
the additive group of $R$ cannot be covered  by finitely many cosets of the
additiv subgroups of $H_{1}, \ldots, H_{t}$;
\item {\rm (ii)}
the group of units of $R$ cannot be  covered by finitely many cosets of the
unit subgroups of $H_{1},\ldots, H_{t}$.
\endproclaim

\proclaim {Theorem 4}
The group of units of a field is a proper union of subsemigroups if and only
if the field is not an algebraic extension of a finite field.
\endproclaim

\proclaim {Remark}
As it was pointed out by I. Ruzsa the polynomial ring $\Bbb Z[x]$ is a proper
union of the following three rings:
$$
\gather
S_1=\{f(x)\in \Bbb Z[x]\mid f(0) \text { is even}\},\quad
  S_2=\{f(x)\in \Bbb Z[x]\mid f(1) \text { is even}\},\\ \vspace{3pt}
S_3=\{f(x)\in \Bbb Z[x]\mid f(0)+f(1) \text { is even}\}.\endgather $$
\endproclaim

\head {3. Preliminaries} \endhead

\proclaim {Lemma 1}
Let  $H_{1}, H_{2}, \ldots, H_{t}$ be subgroups of the group $G$.
If $G$ is covered by finite number of cosets of the $H_{i}$ then at
least one of these subgroups has finite index.
\endproclaim

\demo{Proof}
We use induction on the number of the subgroups.  The statement is
evident if $t=1$ and assume its truth for $t-1$.

We may suppose that $H_{t}$ has infinite index. Then there exists
a coset $H_{t}g$ which
is not in the cover. Hence $H_{t}g$ is covered by finite number of
cosets of $H_{1}, H_{2}, \ldots, H_{t-1}$. If these cosets are
multiplied by $g^{-1}$, a cover of $H_{t}$ is obtained. Thus we
can construct a new cover of $G$ with finite number of cosets of
 $H_{1}, H_{2}, \ldots, H_{t-1}$, and by the
inductive hypothesis Lemma 1 follows.
\hfill $\qed$
\enddemo

\proclaim {Lemma 2 {\rm (\cite{6}, \cite{7})}}
A group is the additive group of an  Artinian ring if and only if it
has the form
$$
\tsize
\bigoplus\limits_{\frak M}{\Bbb Q}\oplus
 \bigoplus\limits_{\text{finite}}C_{p_i^{\infty}}
 \oplus\bigoplus\limits_{\frak N}C_{q_j^{k_j}}, $$
where $p_i$, $q_i$ are prime numbers, ${\frak N}$,
and ${\frak M}$ are
arbitrary cardinals and the factors $q_j^{k_j}$
are divisors of  a fixed natural number $m$.
\endproclaim

\proclaim {Lemma 3 {\rm (\cite{8}, \cite{9})}}
A group is the additive group of a ring with minimal condition for
principal left ideals if and only if its  additive group is
a direct sum of a divisible group and a torsion  group.
\endproclaim

\proclaim {Lemma 4}
Let $\{G_\gamma\mid \gamma\in \Gamma\}$ be a family of abelian groups.
If $G_\gamma$ is not a proper union of finitely many cosets for every
$\gamma$, then $G=\oplus_{\gamma\in\Gamma}G_{\gamma}$ is also not a proper
union of finitely many cosets.
\endproclaim

\demo{Proof}
To prove it by transfinite induction we have two cases to distinguish. If
$\Gamma$ is not a limit ordinal, that is, $\Gamma=\Gamma'+1$ with some
$\Gamma'$ and  for $\Gamma'$ the statement is true. Then we get
$G=G_{\gamma}\oplus G'$, where $G'=\bigoplus_{\gamma'\in \Gamma'}G_{\gamma'}$.
Let $S$ be a coset of $G$ with respect to a subgroup $H$ such that
$b+G_{\gamma}\subseteq S$ with some $b\in G'$. Then $G_{\gamma}\subset H$ and
$S$ has the form $S=G_{\gamma}+S'$, where $S'$ is a proper coset of $G'$.

Suppose that $G$ is a proper union of the cosets $S_1,\cdots,S_n$. If $S_l$
contains a coset of the form $b+G_{\gamma}$ then it can be written as
$G_{\gamma}+S_l'$; otherwise, $S_l'$ is the empty set.  By induction
$$
\bigcup\limits_{l=1}^{n}S_l'\ne G', $$
therefore, there is a $d\in G'$, such that $d+G_{\gamma}$ is not contained in
any $S_l'$. Moreover, if $(d+G_{\gamma})\cap S_l$ is not empty then it
contains an $r_l+d$ and $S_l=r_l+d+G_l$, where $r_l\in G_{\gamma}$ and $G_l$
is a subgroup of $G$. The relations
$$
\gather
S_l\cap (d+G_{\gamma})=(r_l+d+G_l)\cap(r_l+d+ G_{\gamma})
  =(r_l+d)+G_l\cap G_{\gamma}\\ \vspace{3pt}
\tag "and"\\
d+G_\gamma={
\bigcup\limits_{l=1}^{n}}S_l\cap(d+G_{\gamma}) \endgather $$
imply that $G_{\gamma}$ is a proper union of some of the cosets
$r_l+(G_l\cap G_{\gamma})$, which  contradicts.

In the second case $\Gamma$ is a limit ordinal. For a $\Gamma'<\Gamma$  set
$$
G_{\Gamma'}=
\tsize
\bigoplus\limits_{\alpha\in \Gamma'}G_{\alpha}. $$

Assuming $G$ is a proper union of the cosets $T_1, \cdots, T_k$ we obtain
$$
G_{\Gamma'}={
\bigcup\limits_{l=1}^k}(G_{\Gamma}\cap T_l). $$
Since $G_{\Gamma'}\cap T_l$ is also a coset in $G_\gamma$, this
union cannot  be a proper one, that is, for every
$\Gamma'<\Gamma$, $G_{\Gamma'}$ belongs to one of the cosets $T_l$,
$1\leq l\leq k$, which is obviously impossible.
\hfill$\qed$
\enddemo

\demo{Proof {\rm of Theorem 1}}
Let $A$ be the additive group of an Artinian ring. According to Lemma 2, if
the non-divisible part $\oplus_{\frak N}C_{p_i^{k_i}}$ of $A$ contains a
direct summand $C_{p^k_i}$ at least twice, then
$A=L\oplus C_{p^k}\oplus C_{p^l}$ and if $k\leq l$,
$C_{p^k}=\langle a\rangle$, $C_{p^l}=\{b_1, b_2, \ldots, b_{p^l}\}$  and
$A_i=\gp{ L, ab_i}$, therefore,  $A$ is a proper union of the subgroups
$A_1, A_2, \ldots, A_{p^l}$. Furthermore, $\oplus_{\frak N}C_{p_i^{k_i}}$ is a
finite cyclic group.

A quasycyclic group and the additive group of $\Bbb Q$ have no
maximal subgroups, hence by Lemma 1 they are not a proper union of
cosets.  Applying Lemma  2 we may assume that  the maximal divisible subgroup
$B$ of $A$ is not a proper union of cosets.  Clearly, the finite cyclic
group  $C$ is not a proper union of subgroups. It yields that
$A=B\oplus C$ is also not a proper union of subgroups.
\hfill $\qed$
\enddemo

One can repeat the argument detailed above to prove Theorem 2.

\demo{Proof {\rm of Theorem 3}}
(i) If  $R$  is  covered by finitely many  cosets of the additive subgroups
$H_{1}, H_{2}, \ldots, H_{t}$ then by Lemma 1 there exists a subgroup $H=H_i$
of finite index in the additive group of $R$. Let
$$
R=a_1+H\cup a_2+H\cup \ldots\cup a_s+H  $$
be a decomposition of $R$ with respect to $H$. Then $H$ is an infinite subskew
field and in the  infinite set $\{a_i+a_j\lambda\mid 0\ne \lambda\in H \}$
there exist two different $a_i+a_j\lambda_1$ and $a_i+a_j\lambda_2$, which
belong to the same coset $a_k+H$.  Then $a_j(\lambda_1-\lambda_2)\in H$
and $a_j\in H$, which is impossible.

\smallskip
(ii)
Let the group of units $U(R)$  be  covered by finitely many cosets of the
multiplicative subgroups $U(H_{1}),U(H_{2}),\ldots,U(H_{t})$. Then by Lemma~1
there exists a subgroup $H=U(H_i)$ of finite index in the group of units
$U(R)$ and we have the decomposition
$$
U(R)=a_1H\cup a_2H\cup \ldots\cup a_s H.  $$
Then $H$ is an infinite subskew field and in the infinite set
$\{a_i+a_j\lambda\mid 0\ne \lambda\in H \}$ there exist two different
$a_i+a_j\lambda_1$ and $a_i+a_j\lambda_2$, which belong to the same coset
$a_kH$. Therefore, $a_i+a_j\lambda_2=(a_i+a_j\lambda_1)\lambda_3$ and we obtain
$$
a_i(1-\lambda_3)=a_j(\lambda_1\lambda_3-\lambda_2),  $$
and $1-\lambda_3, \lambda_1\lambda_3-\lambda_2\in H$, which is impossible.
\hfill $\qed$
\enddemo
\demo{Proof {\rm of Theorem 4}}
Let $F$ be a field with multiplicative group $U(F)$. If $F$ is an algebraic
extension of a finite field $F_0$ and $U(F)$ is a proper union of the
subsemigroups $M_1, \ldots, M_n$, then there are  elements $m_i\in M_i$ with
$$
m_i\notin\bigcup\limits_{l\ne i}M_l, $$
where $i=\!1, \ldots, n$
furthermore, the multiplicative group of $F_0(m_1,\ldots, m_n)\!$
is a proper union of the groups $M_i\cap F_0(m_1, \ldots, m_n)$,
($i=1, \ldots, n$). However, 
$F_0(m_1,\ldots, m_n)$ is a finite field having
cyclic multiplicative group, which cannot be a proper union.

If $F$ is not an algebraic extension of a finite field then $U(F)$ contains
two multiplicatively independent elements denoted by $z_1$ and $z_2$. Indeed,
if $\text {char}(F)=0$ then one can take $z_1=2$ and $z_2=3$, say; and if $F$
has a transcendental element $\tau$ (over a finite ground field contained in
$F$), then put $z_1=\tau$ and $z_2=\tau+1$. Let $G$ be a multiplicatively
independent generating set for $U(F)$ containing $z_1$ and $z_2$. Moreover,
for a $z\in U(F)$ let $e_i(z)$ ($i=1, 2$) denote the exponent of $z_i$
($i=1,2$) in the expression of~$z$ as a product of generators from $G$.
The lattice $\Bbb Z\oplus \Bbb Z$ is a proper union of the lattices
$L_1, L_2$ and $L_3$ spanned by
$$
\{ (1, 0), (1, 2)\},\, \{(0, 1), (2, 1)\},\, \{(1, 1), (-1, 1)\} $$
respectively, hence $U(F)$ is a proper union of the subsemigroups
$$
\{z\in U(F) \mid (e_1(z), e_2(z))\in L_i\},  $$
where $i=1, 2, 3$.
\hfill $\qed$
\enddemo


\newpage

\Refs



\ref\no{1}
\by K. Venkatachaliengar,   T. Soundararajan
\paper  A covering theorem for skew fields
\jour Indag. Math.
\vol 31
\yr 1969
\pages 441--442
\endref

\ref\no{2}
\by K. Kishimoto,   K. Motose
\paper A covering theorem for simple rings
\jour J. Fac. Sci. Shinshu Univ.
\vol 5
\yr 1970
\pages 107--108
\endref

\ref\no{3}
\by M. \^Ohori
\paper On finite unions of subrings
\jour  Math. J. Okayama Univ.
\vol 19
\yr 1976/77
\pages 47--50
\endref

\ref\no{4}
\by T. Laffey
\book Review
\bookinfo  MR:10512
\yr 1978
\endref

\ref\no{5}
\by J. Lewin
\paper Subrings of finite index of finitely  generated rings
\jour J. Algebra
\vol 5
\yr 1967
\pages 84--88
\endref

\ref\no{6}
\by  L. Fuchs,   T. Szele
\paper On Artinian rings
\jour Acta Sci. Math. Szeged
\vol 17
\yr 1956
\pages 30--40
\endref

\ref\no{7}
\by  A. Kert\'esz
\book Lectures on Artinian rings
\publ Budapest
\publaddr Akad\'emiai Kiad\'o
\yr 1987
\endref

\ref\no{8}
\by F. A. Sz\'asz
\paper \"Uber Artinsche Ringe
\jour Bull. Acad. Polon. Sci.
\vol  11
\yr 1963
\pages 351--354
\endref


\ref\no{9}
\by  L. Fuchs
\book Abelian groups, II
\publ Academic Press
\publaddr New York and London
\yr 1973
\endref

\endRefs
\enddocument